\documentclass[12pt]{article}
\usepackage{color}
\usepackage{amsmath}
\usepackage{amssymb}
\usepackage{amsbsy}

\usepackage{amsfonts}
\usepackage{tikz}

\definecolor{darkblue}{rgb}{0,0,0.8}
\definecolor{alizarin}{rgb}{0.82, 0.1, 0.26}

\newtheorem{lem}{Lemma}[section]%
\newtheorem{thm}[lem]{Theorem}%
\newtheorem{defi}[lem]{Definition}%
\newtheorem{ques}{Question}%
\newtheorem{prop}[lem]{Proposition}%

\def\a{\alpha} \def\b{\beta}   
 \def\s{\sigma}

\def\G{\Gamma}

 \def\lg{\langle} \def\rg{\rangle}
\def\nd{\mathrel{\bigm|\kern-.7em/}}

\def\f{\noindent}

\def\PSL{\hbox{\rm PSL}}
\def\SL{\hbox{\rm SL}} \def\GL{\hbox{\rm GL}}

\DeclareMathOperator{\Aut}{Aut}

\def\soc{\hbox{\rm soc}}

\def\Cay{\hbox{\rm Cay}}

\def\mz{{\mathbb Z}}
\def\K{{\rm\bf K}}

\def\H{\mathcal{H}}
\def\P{\mathcal{P}}

\def\R{\mathcal{R}}
\def\demo{\noindent{\bf Proof}\hskip10pt}

\def\qed{\hskip10pt $\Box$\vspace{3mm}}

\tikzset{every picture/.style={line width=0.75pt}} 

\begin{document}
\title{Using mixed dihedral groups to construct normal Cayley graphs, and a new bipartite $2$-arc-transitive graph which is not a Cayley graph}
\author
{
Daniel R. Hawtin\\
{\small Faculty of Mathematics, The University of Rijeka}\\[-4pt]
{\small Rijeka, 51000, Croatia}\\[-2pt]
{\small{Email}:\ dan.hawtin@gmail.com}\\[+6pt]
Cheryl E. Praeger\\
{\small Department of Mathematics and Statistics, The University of Western Australia}\\[-4pt]
{\small Crawley, WA 6907, Australia}\\[-2pt]
{\small{Email}:\ cheryl.praeger@uwa.edu.au}\\[+6pt]
Jin-Xin Zhou  \\
{\small Department of Mathematics, Beijing Jiaotong University}\\[-4pt]
{\small Beijing 100044, P.R. China}\\[-2pt]
{\small{Email}:\ jxzhou@bjtu.edu.cn}
}

\date{}
\maketitle
\begin{abstract}
A \emph{mixed dihedral group} is a group $H$ with two disjoint subgroups $X$ and $Y$, each elementary abelian of order $2^n$, such that $H$ is generated by $X\cup Y$, and $H/H'\cong X\times Y$. In this paper we give a sufficient condition such that the automorphism group of the Cayley graph $\Cay(H,(X\cup Y)\setminus\{1\})$ is equal to $H: A(H,X,Y)$, where $A(H,X,Y)$ is the setwise stabiliser in $\Aut(H)$ of $X\cup Y$. We use this criterion to resolve a questions of Li, Ma and Pan from 2009, by constructing a $2$-arc transitive normal cover of order $2^{53}$ of the complete bipartite graph $\K_{16,16}$ and prove that it is \emph{not} a Cayley graph.

\end{abstract}

\section{Introduction}
In \cite{HPZ}, the authors introduced the family of \emph{mixed dihedral groups} and used them to construct new families of graphs with a high degree of symmetry, shedding light on several questions raised in the literature. Each mixed dihedral group $H$ has two specified elementary abelian $2$-subgroups $X, Y$ such that $X\cup Y$ generates $H$, the derived quotient $H/H'\cong X\times Y\cong C_2^n\times C_2^n$, and studying links between the Cayley graph $\Gamma(H,X,Y)$ for $H$ relative to $X, Y$, and its clique graph $\Sigma(H,X,Y)$, has proved interesting and fruitful. For example, one family of mixed dihedral groups yielded an infinite family of normal Cayley graphs $\Gamma(H,X,Y)$ which were $2$-geodesic transitive, but neither distance transitive nor $2$-arc-transitive (\cite[Theorem 1.8]{HPZ}), while a different family of these groups produced semisymmetric, locally $2$-arc-transitive normal covers $\Sigma(H,X,Y)$ of complete bipartite graphs (that is, these graphs were not vertex-transitive, \cite[Theorem 1.3]{HPZ-paper2}). In each case the graphs provided examples of highly symmetric graphs with properties not seen previously: for example, the semisymmetric graphs were the first such graphs constructed from groups which were not $2$-generated.

In this paper we obtain a sufficient condition for the Cayley graphs $\Gamma(H,X,Y)$ of mixed dihedral $2$-groups  to be normal (Theorem~\ref{mix-dihedrant-characterization}) shedding light on a question asked in 1999 by the second author \cite[Question 4]{Praeger-NE-Cay} about the automorphism groups of Cayley graphs. We then examine the clique graph $\Sigma(H,X,Y)$ of order $2^{53}$ for a particular mixed dihedral group $H$, and use Theorem~\ref{mix-dihedrant-characterization} to prove that it is a $2$-arc-transitive normal cover of $\K_{16,16}$,  but not a Cayley graph (Theorem~\ref{non-Cayley}).
This answers in the negative a question asked in 2009 by C. H.~Li.

\begin{ques}{\rm\cite[Problem~(1)]{LMP2009}}\label{prob-1}
Is every vertex-transitive and locally primitive normal cover of $\K_{2^n,2^n}$ a Cayley graph?
\end{ques}

\subsection{Cayley graphs}\label{s:Cay}

Let $G$ be a finite group $G$ and  $S\subseteq G\setminus\{1\}$. Then $S$ is \emph{inverse-closed} if $s^{-1}\in S$ for all $s\in S$. The {\em Cayley graph} $\G:=\Cay(G,S)$ of $G$ with respect to an inverse-closed subset $S$ of $G\setminus\{1\}$ is the graph with vertex set $G$ and edge set $\{\{g,sg\} \mid g\in G,s\in S\}$. This graph is connected if and only if $S$ generates $G$. For any $g\in G$ define
\[
\hat{g}: x\mapsto xg\ \mbox{for $x\in G$ and set $\hat{G}=\{\hat{g}\ \mid g\in G\}$.}
\]
Then $\hat{G}$ is a regular permutation group on $V(\G)$ (that is, $\hat{G}$ is transitive and all vertex-stabilisers are trivial, see for example \cite[Lemma 3.7]{PS}) and $\hat{G}$ is a subgroup of $\Aut(\G)$, and it is convenient to identify $\hat{G}$ with $G$. In fact, a graph is a Cayley graph if and only if its automorphism group has a subgroup that is regular on the vertex-set (see \cite[Lemma~16.3]{Biggs}).
Let
\begin{equation}\label{eq-ags}
    \Aut(G,S)=\{\a\in\Aut(G): S^\a=S\}.
\end{equation}

Then the normaliser of $G$ in $\Aut(\G)$ is $G\rtimes \Aut(G,S)$ (see~\cite{Godsil-1981}). A Cayley graph $\G=\Cay(G,S)$ is said to be {\em normal} if  $\Aut(\G)$ is equal to $G\rtimes \Aut(G,S)$ (see \cite{Xu98}), and to be  {\em normal-edge-transitive} if $G\rtimes \Aut(G,S)$ is transitive on $E(\G)$ (see \cite{Praeger-NE-Cay}). In particular for a normal-edge-transitive Cayley graph $\G=\Cay(G,S)$, \cite[Question 4]{Praeger-NE-Cay} asked for conditions on $G, S$ under which we could `guarantee that $G\rtimes \Aut(G, S)$ is not much smaller than $\Aut(\G)$'.  Theorem~\ref{mix-dihedrant-characterization} contributes to our understanding of these graphs.




\subsection{Mixed dihedral groups and their graphs}\label{s:mix-dih}

These groups have an associated integer parameter $n\geq2$.

\medskip\noindent
{\rm (a)} A finite group $H$ is an {\em $n$-dimensional mixed dihedral group relative to $X$ and $Y$} if $X$ and $Y$ are subgroups of $H$ such that $X\cong Y\cong C_2^n$, $H=\lg X, Y\rg$, and $H/H'\cong C_2^{2n}$, where $H'$ is the derived subgroup of $H$.

    \medskip\noindent
{\rm (b)} For $H, X, Y$ as in part (a),  the graphs $\Gamma(H,X,Y)$ and $\Sigma(H,X,Y)$ are defined as follows:
\begin{equation}\label{eq-1}
\Gamma=\Gamma(H,X,Y)=\Cay(H,S(X,Y)),\ {\rm with}\ S(X,Y)=(X\cup Y)\setminus\{1\};
\end{equation}
and $\Sigma=\Sigma(H,X,Y)$ is the graph with vertex-set and edge-set given by:
\begin{equation}\label{eq-2}
\begin{array}{l}
V(\Sigma)=\{Xh, Yh: h\in H\},\\
E(\Sigma)=\{\{Xh, Yg\}:   h,g\in H\ \mbox{and}\ Xh\cap Yg\neq\emptyset\}.
\end{array}
\end{equation}
Because of the relationship between the set $S(X,Y)$ and the subgroups $X, Y$, it is convenient to denote
\begin{equation}\label{eq-ahxy}
    \text{$A(H,X,Y):=\Aut(H, S(X,Y))$ as in \eqref{eq-ags}, and let $G:=H\rtimes A(H,X,Y).$}
\end{equation}







Our first main result is a criterion which guarantees that the Cayley graphs $\G(H,X,Y)$ are normal.

\begin{thm}\label{mix-dihedrant-characterization}
Let $n\geq 2$, and let $H$ be an $n$-dimensional mixed dihedral group relative to $X$ and $Y,$ such that $H$ is a $2$-group. Suppose that $\G=\G(H,X,Y)$ is normal edge-transitive, that is, the group $G$ in \eqref{eq-ahxy}  is transitive on $E(\G)$. If the pointwise stabiliser $A(H,X,Y)_{(X)}$ is non-trivial, and $A(H,X,Y)$ does not contain $\Aut(X)\times\Aut(Y)$, then $G=\Aut(\G)$.
\end{thm}

\subsection{$2$-arc-transitive normal covers of complete bipartite graphs}\label{s:sigma}

In \cite{LMP2009}, Li et al. characterised vertex-transitive, locally primitive graphs of prime power order. More precisely, they proved that every such graph is either a Cayley graph or a normal cover of the complete bipartite graph $\K_{2^n,2^n}$. Further, they asked  whether all of the normal covers in the second case were also Cayley graphs (see Question~\ref{prob-1}).

In this second part of the paper we present an example which we prove is not a Cayley graph, thus answering this question in the negative. Our example is quite large and complicated, arising as a graph $\Sigma(H,X,Y)$ for a certain mixed dihedral group $H$. It raises several questions:

\begin{ques}\label{prob-2}
Are there infinitely many vertex-transitive, locally primitive normal covers of some complete bipartite graph $\K_{2^n,2^n}$ which are not Cayley graphs?
What is the smallest such graph?
\end{ques}

The automorphism group of $\Delta=\K_{2^n,2^n}$ is the wreath product $\Aut(\Delta)=S_{2^n}\wr S_2$. For a subgroup $L\leq \Aut(\Delta)$, the graph $\Delta$ is  said to be {\em $L$-edge-affine} if $L\cap (S_{2^n}\times S_{2^n})$ has a normal subgroup $C_{2}^n\times C_{2}^n$ that is regular on $E(\Delta)$. Two theorems we proved in \cite{HPZ} concerning normal covers of $\K_{2^n,2^n}$ provided the clue which ultimately led us to our construction of our non-Cayley example. Namely, we proved:

\begin{quote}
\cite[Theorem~1.4]{HPZ} \quad If $\Sigma$ is a $(G,2)$-arc-transitive graph and is a normal cover of $\K_{2^n,2^n}$, then either $\Sigma$ is a Cayley graph, or $n\geq4$ and the normal quotient $\K_{2^n,2^n}$ is $L$-edge-affine relative to the group $L$ induced by $G$.
\end{quote}
Moreover, in the latter case, we showed in \cite[Theorem~1.6]{HPZ} (see Theorem~\ref{iff-characterization}) that $G$ has a normal subgroup $H$ which is an $n$-dimensional mixed dihedral group relative to $X$ and $Y$, and $G$ is edge-transitive on the line graph of $\Sigma$, which is isomorphic to $\G(H,X,Y)$. Thus we needed to search for constructions involving mixed dihedral $2$-groups. Our graph is based on the following $2$-group.




\begin{defi}\label{def:H}
Let $X_0=\{x_1,x_2,x_3,x_4\}$, $Y_0=\{y_1,y_2,y_3,y_4\}$, and let
\[
\P=\lg X_0\cup Y_0\cup\{r\}\mid \mathcal{R}\rg
\]
where $\mathcal{R}$ is the following set of relations: for $x,x'\in X_0$, $y,y'\in Y_0$,  $z,z', z'', z''' \in X_0\cup Y_0$ and $\sigma=(1,2,4,3)$,
\[\begin{array}{l}
r^8=1, z^2=1, [x,x']=[y,y']=1,[x,y]^2=1, [[x,y],z]^2=1, [[[z,z'],z''],z''']=1,\\

x_i^r=y_i, \ \text{and}\ y_i^r=x_{i\sigma}, \ \text{for $1\leq i\leq 4$, }  \\

[[x_1x_2x_3x_4,y_1],x_1]\cdot[[y_1,x_1],y_2]\cdot[[x_1,y_2],x_3]\cdot[[y_2,x_3],y_4]\cdot[[x_3,y_4],x_2]\cdot[[y_4,x_2],y_3]=\\

[[x_2,y_3],x_1x_2x_3x_4]\cdot[[y_3,x_1x_2x_3x_4],y_1],\\

[[x_4,y_2], x_1x_2x_3x_4]\cdot[[y_2,x_1x_2x_3x_4],y_3]\cdot[[x_1x_2x_3x_4,y_3],x_2]\cdot[[y_3,x_2],y_1y_2y_3y_4]=\\

[[x_2,y_1y_2y_3y_4],x_1]=[[y_1y_2y_3y_4,x_1],y_4]\cdot[[x_1,y_4],x_4]\cdot[[y_4,x_4],y_2].
\end{array}\]
\end{defi}

Here is our result, answering Question~\ref{prob-1} in the negative.

\begin{thm}\label{non-Cayley}
Let $\P, X_0, Y_0$ be as in Definition~\emph{\ref{def:H}} and let  $X=\lg X_0\rg$, $Y=\lg Y_0\rg$, and  $\H=\lg X_0\cup Y_0\rg$.  Then the following hold.
\begin{enumerate}
  \item [{\rm (1)}]\ $\P= \H\rtimes\lg r\rg=\H\rtimes C_8$, $|\P|=2^{59}$, and $\H$ is a $4$-dimensional mixed dihedral group of order $2^{56}$ relative to $X$ and $Y$.

  \item [{\rm (2)}]\ $\G=\G(\H,X,Y)$ as in \eqref{eq-1} is a normal Cayley graph and is edge-transitive, with
 $\Aut(\G)=\H\rtimes A(\H,X,Y)$ as in \eqref{eq-ahxy} and
    $A(\H,X,Y)\cong (C_{15}\times C_{15}): C_8$.

  \item [{\rm (3)}]\ $\Sigma=\Sigma(\H, X, Y)$ is a $2$-arc-transitive normal cover, of order $2^{53}$, of $\K_{2^4,2^4}$, and $\Sigma$ is not a Cayley graph.
\end{enumerate}
\end{thm}

\section{Preliminaries}\label{prelimsect}

All graphs we consider are finite, connected, simple and undirected.
Let $\G$ be a graph. Denote by $V(\Gamma)$, $E(\Gamma)$ and $\Aut(\Gamma)$ the vertex set, edge set, and full automorphism group of $\G$, respectively. For $v\in V(\G)$, let $\G(v)$ denote the set of vertices adjacent to $v$. A graph $\G$ is said to be {\em regular} if there exists an integer $k$ such that $|\G(v)|=k$ for all vertices $v\in V(\G)$. A graph $\G$ is bipartite if $E(\G)\ne\emptyset$ and $V(\G)$ is of the form $\Delta\cup\Delta'$ such that each edge consists of one vertex from $\Delta$ and one vertex from $\Delta'$. If $\G$ is connected then this vertex partition is uniquely determined and the two parts $\Delta, \Delta'$ are often called the \emph{biparts} of $\G$.

Let $G\leq\Aut(\G)$. For $v\in V(\G)$, let $G_v=\{g\in G\ :\ v^g=v\}$, the stabiliser of $v$ in $G$. We say that $\G$ is {\em $G$-vertex-transitive\/} or {\em $G$-edge-transitive\/} if $G$ is transitive on $V(\G)$ or $E(\G)$, respectively.
When $G=\Aut(\G)$, a $G$-vertex-transitive or $G$-edge-transitive graph $\G$ is simply called {\em vertex-transitive\/} or {\em edge-transitive\/}, respectively. A regular graph $\G$ is said to be {\em $G$-locally primitive\/} or {\em locally $(G,2)$-arc-transitive\/} if $G\leq\Aut(\G)$ and $G_v$ is primitive or $2$-transitive  on $\G(v)$, respectively, for each $v\in V(\G)$. Similarly, when $G=\Aut(\G)$, a $G$-locally primitive or locally $(G,2)$-arc-transitive graph $\G$ is simply called {\em locally primitive\/} or {\em locally $2$-arc-transitive\/}, respectively.

\subsection{Normal quotients and normal covers of graphs}\label{s:nquots}
 The normal quotient method for investigating vertex- or edge-transitive graphs proceeds as follows.
 Assume that $G\leq\Aut(\G)$ is such that $\G$ is $G$-vertex-transitive or $G$-edge-transitive. Let $N$ be a normal subgroup of $G$ such that $N$ is intransitive on $V(\G)$. The {\em $N$-normal quotient graph\/} of $\G$ is defined as the graph $\G_N$ with vertices the $N$-orbits in $V(\G)$ and with two distinct $N$-orbits adjacent if there exists an edge in $\G$ consisting of one vertex from each of these orbits. If $\G$ is regular, and if $\G_N$ and $\G$ have the same valency, then we say that $\G$ is an {\em $N$-normal cover} of $\G_N$.

The following lemma is about the normal quotients of locally primitive graphs.

\begin{lem}{\rm\cite[Lemma~2.2]{HPZ}}\label{quot}
Let $\G$ be a connected $G$-locally primitive bipartite graph of valency $k>1$ with  bipartition $V(\G)=O_1\cup O_2$, so each $|O_i|>1$. Suppose that $N\unlhd G$ is such that $N$ fixes both $O_1$ and $O_2$ setwise, and $N$ is intransitive on $O_1$ and on $O_2$. Then
\begin{enumerate}
  \item [{\rm (1)}]\ $N$ acts semiregularly on $V(\G)$, $N$ is the kernel of $G$ acting on $V(\G_N)$ and $G/N\leq\Aut(\G_N)$.

  \item [{\rm (2)}]\ For $N\leq H\leq G$, $H$ is regular on $E(\G)$  if and only if $H/N$ is regular on $E(\G_N)$.
\end{enumerate}
\end{lem}

\subsection{Cliques, clique graphs and line graphs}

A \emph{clique} of a graph $\G$ is a subset  $U\subseteq
V(\G)$ such that every pair of vertices in $U$ forms an
edge of $\G$. A clique $U$ is \emph{maximal} if no subset
of $V(\G)$  properly containing $U$ is a clique. The \emph{clique graph}  of $\G$ is defined as the graph $\Sigma(\G)$ with vertices the maximal cliques of $\G$ such that two distinct maximal cliques are adjacent in $\Sigma(\G)$ if and only if their intersection is non-empty. Similarly the \emph{line graph}  of $\G$ is defined as the graph $\mathcal{L}(\G)$ with vertex set $E(\G)$ such that two distinct edges $e,e'\in E(\G)$ are adjacent in $\mathcal{L}(\G)$ if and only if $e\cap e'\ne\emptyset$.

\subsection{Cayley graphs and bi-Cayley graphs}\label{sub:cay}

Recall that a group $G$ of permutations of a set $V(\G)$ is called \emph{regular} if it is transitive, and some (and hence all) stabilisers $G_v$ are trivial. More generally $G$ is called \emph{semiregular} if the stabiliser $G_v=1$ for all $v\in V(\G)$. So $G$ is regular if and only if it is semiregular and transitive.

As discussed in Section~\ref{s:Cay}, for a Cayley graph $\G=\Cay(G,S)$, the group $\widehat{G}$ of right translations by elements of $G$ is regular on $V(\G)$, and
Cayley graphs are precisely those graphs $\Gamma$ for which $\Aut(\G)$ has a subgroup that is regular on $V(\G)$ (see \cite[Lemma~16.3]{Biggs}). Similarly a  graph $\Gamma$ is called a \emph{bi-Cayley graph} if $\Aut(\G)$ has a subgroup $H$ which is semiregular on $V(\G)$ with two orbits. The following results from \cite{HPZ} will be useful.

\begin{lem}{\rm\cite[Lemma~2.3]{HPZ}}\label{bi-Cayley}
Let $\G$ be a connected bi-Cayley graph of a group $H$ such that the two orbits of $H$ on $V(\G)$ contain no edges of $\G$ and let $N=N_{\Aut(\G)}(H)$. Then, for each $v\in V(\Gamma)$, the stabiliser $N_v$ acts faithfully on $\Gamma(v)$.
\end{lem}

\begin{lem}{\rm\cite[Lemma~2.6]{HPZ}}\label{sufficient}
Let $\G$ be a connected $(G,2)$-arc-transitive graph, and let $u\in V(\G)$. Suppose that $\G$ is an $N$-normal cover of $\K_{2^n,2^n}$, for some normal $2$-subgroup $N$ of $G$, and that the stabiliser $G_u$ acts non-faithfully on $\G(u)$. Then $N\unlhd\Aut(\G)$.
\end{lem}

\subsection{Mixed dihedrants and their clique graphs}

We record the properties we will need of the graphs defined in Subsection~\ref{s:mix-dih}.

\begin{lem}{\rm\cite[Lemmas~4.1-4.2]{HPZ}}\label{lem:prop-mixed-dih}
Let $H$ be an $n$-dimensional mixed dihedral group relative to $X$ and $Y$ with $|X|=|Y|=2^n\geq4$, let $\G=\G(H,X,Y)$ and $\Sigma=\Sigma(H,X,Y)$ be the graphs  defined in \eqref{eq-1} and~\eqref{eq-2}, and let  $G=H:A(H,X,Y)$ as in \eqref{eq-ahxy}. Then the following hold.
\begin{enumerate}
\item [{\rm (1)}] $\Sigma$ is the clique graph of $\G$, and the  map $\varphi:z\to \{Xz,Yz\}$, for $z\in H$, induces an isomorphism from  $\G$ to the line graph $\mathcal{L}(\Sigma)$ of $\Sigma$.


\item [{\rm (2)}] $\Aut(\G)=\Aut(\Sigma)= \Aut(\mathcal{L}(\Sigma))$.

 \item [{\rm (3)}] $G\leq \Aut(\Sigma)$ is edge-transitive;  and  $H$ is regular on $E(\Sigma)$ with two orbits on $V(\Sigma)$.


\item [{\rm (4)}] $A(H, X, Y))\cong A(H,X,Y)^{X\cup Y}\leq (\Aut(X)\times\Aut(Y)): C_2\cong(\GL(n,2)\times\GL(n,2)): C_2$ where the $C_2$ interchanges $X$ and $Y$.

\end{enumerate}
\end{lem}

\begin{thm}\label{iff-characterization}{\rm\cite[Theorem~1.6]{HPZ}}
Let $n\geq 2$, and let $\Sigma$ be a graph,  $G\leq \Aut(\Sigma)$, and $N\lhd G$. Then the following are equivalent.
 \begin{enumerate}
    \item[{\rm (a)}]  $\Sigma$ is a $(G,2)$-arc-transitive $N$-normal cover of $\K_{2^n,2^n}$ which is $G/N$-edge-affine;
    \item[{\rm (b)}] $G$ has a normal subgroup $H$ with $H'=N$ such that $H$ is an $n$-dimensional mixed dihedral group relative to $X$ and $Y$, the line graph of $\Sigma$ is $\G(H,X,Y)$ as in \eqref{eq-1}, and $\G(H,X,Y)$ is $G$-edge-transitive.
 \end{enumerate}
\end{thm}

\subsection{$2$-arc-transitive subgroups of $\Aut(\K_{2^n,2^n})$}

The following proposition will be used in the proof of Theorem~\ref{mix-dihedrant-characterization}.

\begin{prop}{\rm\cite[Proposition~3.1]{HPZ}}
\label{basic-auto}
Let $n\geq 2$ be a positive integer and let $\G=\K_{2^n,2^n}$ with biparts $U$ and $W$. Let $G\leq\Aut(\G)$ be $2$-arc-transitive on $\G$, and let $G^+=G_U=G_W$, the setwise stabiliser in $G$ of each of the biparts. Then one of the following holds.
\begin{enumerate}
  \item [{\rm(1)}]\ $\soc(G^+)\cong A_{2^n}\times A_{2^n}$ with $n\geq3$;
  \item [{\rm(2)}]\ $\soc(G^+)\cong\PSL(2,p)\times\PSL(2,p)$ with $p=2^n-1$ a prime and $n\geq 3$;
  \item [{\rm(3)}]\ $\soc(G^+)\cong C_2^{2n}$,  acting regularly on the edge set of $\G$;
  \item [{\rm(4)}]\ $n=3$ and $G^+={\rm AGL(3,2)}, G_u=\PSL(2,7), G_w=\SL(3,2)$, with $G_{u}\cap G_w=\mz_7:\mz_3$, where $u\in U$ and $w\in W$.
\end{enumerate}
Moreover, if $G$ does not contain a  subgroup acting regularly on $V(\G)$, then $n\geq 4$ and
case {\rm(3)} holds.
\end{prop}

\f\section{Proof of Theorem~\ref{mix-dihedrant-characterization}}

Suppose that $n\geq 2$, and that $H$ is an $n$-dimensional mixed dihedral group relative to $X$ and $Y,$ such that $H$ is a $2$-group. Suppose also that
$G:=H\rtimes A(H,X,Y)$, as in \eqref{eq-ahxy}, is edge-transitive on $\G=\G(H,X,Y)$, and that the pointwise stabiliser $A(H,X,Y)_{(X)}$ is non-trivial, and $A(H,X,Y)$ does not contain $\Aut(X)\times\Aut(Y)$. We need to prove that $G=\Aut(\G)$, but we assume to the contrary that $G<\Aut(\G)$, and derive a contradiction.

By Lemma~\ref{lem:prop-mixed-dih}~(1), $\G$ is the line graph of the graph $\Sigma=\Sigma(H,X,Y)$, as in \eqref{eq-2}. Thus the group $G$ satisfies all the conditions of Theorem~\ref{iff-characterization}(b) and hence, by that result, $\Sigma$ is a $(G,2)$-arc-transitive $H'$-normal cover of $\K_{2^n,2^n}$, and the quotient $\Sigma=\Sigma(H,X,Y)$ is $G/H'$-edge-affine (as defined in Subsection~\ref{s:sigma}). In addition, $H'$ is semiregular on $V(\Sigma)$, by Lemma~\ref{quot}(1). Since by assumption $A(H,X,Y)_{(X)}\ne 1$, the stabiliser $G_X$ of the vertex $X\in V(\Sigma)$ is not faithful on the neighbourhood $\{Yx : x\in X\}$ of $X$ in $\Sigma$. Since $H'$ is a normal $2$-subgroup of $G$, it follows from Lemma~\ref{sufficient} that  $H'\unlhd\Aut(\Sigma)$. Thus by \cite[Theorem 4.1]{ONS} applied to $G$ and $\Aut(\Sigma)$, each of $G/H'$ and $\Aut(\Sigma)/H'$ is $2$-arc-transitive on the normal quotient $\Sigma_{H'}\cong \K_{2^n,2^n}$. Now we apply Proposition~\ref{basic-auto} to $\Sigma_{H'}$ and $\Aut(\Sigma)/H'$. For any subgroup $T$ satisfying $H'\leq T\leq  \Aut(\Sigma)$, let $T^+$ denote the subgroup of $T$ such that $T^+/H'$ is the setwise stabiliser in $T/H'$ of each of the parts of the bipartition of $\Sigma_{H'}$. By Lemma~\ref{lem:prop-mixed-dih}~(3), $H=H^+\leq \Aut(\Sigma)^+$, and $H$ is regular on $E(\Sigma)$, whence by  Lemma~\ref{quot}~(2),  $H/H'$ is regular on $E(\Sigma_{H'})$.  Recall that $H/H'\cong C_2^{2n}$ by the definition of a mixed dihedral group $H$.

\smallskip\noindent
{\em Proposition~\ref{basic-auto} cases~(2) or (4):}\quad Here the group $\Aut(\Sigma)^+/H'$ does not contain a subgroup isomorphic to $C_2^{2n}$, so these cases do not arise.

\smallskip\noindent
{\em Proposition~\ref{basic-auto} case~(1):}\quad
Here $\soc(\Aut(\Sigma)^+/H')\cong A_{2^n}\times A_{2^n}$ with $n\geq3$. Since $H/H'$ is regular on the edge set of $\Sigma_N\cong \K_{2^n,2^n}$, we have $H/H'=(H/H')^{U}\times (H/H')^{W}$, where $U,W$ are the two parts of the bipartition of $\Sigma_{H'}$, and $(H/H')^U$ and $(H/H')^{W}$ are the permutation groups induced by $H/H'$ on $U$ and $W$, respectively. Then $(H/H')^U$ acts regularly on $U$ and $(H/H')^W$ acts regularly on $W$. The permutation induced on $U$ by a nontrivial element of $(H/H')^U$ is a product of $2^{n-1}$ cycles of length $2$, and hence is an even permutation of $U$ since $n\geq3$. Similarly every element of $(H/H')^W$ induces an even permutation of $W$.
It follows that $H/H'\leq \soc(\Aut(\Sigma)^+/H')\cong A_{2^n}\times A_{2^n}$, and so the normaliser of $H/H'$ in $\soc(\Aut(\Sigma)^+/H')$ is isomorphic to $H/H': (\GL(n,2)\times \GL(n,2))$. This implies that $\Aut(X)\times\Aut(Y)\leq A(H,X,Y)$, contradicting our assumption.


\smallskip\noindent
{\em Proposition~\ref{basic-auto} case~(3):}\quad In this final case $T:=\soc(\Aut(\Sigma))^+$ satisfies $T/H'\cong C_2^{2n}$, and $T/H'$ acts regularly on  $E(\Sigma_{H'})$.

\smallskip\noindent
{\em Claim: $H/H'$ is a minimal normal subgroup of $G/H'$.}\quad Let $H'< M \leq H$ such that $M/H'$ is a minimal normal subgroup of $G/H'$ (so $M\unlhd G$ and $M/H'\leq H/H'$). As shown above, $G/H'$ is $2$-arc-transitive on $\Sigma_{H'}\cong\K_{2^n,2^n}$, and hence $G^+/H'$ acts $2$-transitively on each of  the two parts  $U,W$ of the bipartition of $\Sigma_{H'}$. Thus $M/H'$ is transitive on each of $U$ and $W$. As noted above, $H'$ is semiregular on $V(\Sigma)$. Hence, if $M/H'$ is semiregular on $V(\Sigma_{H'})$, then $M$ is  semiregular on $V(\Sigma)$ and its orbits are the two parts of the bipartition of $V(\Sigma)$. This implies that $\Sigma$ is a bi-Cayley graph of $M$ (as defined in Section~\ref{sub:cay}), and since $M\unlhd G$ it follows from Lemma~\ref{bi-Cayley} that $G_X$ acts faithfully on the neighborhood of $X$ in $\Sigma$, which is a contradiction. Thus $M/H'$ is not semiregular on $V(\Sigma_{H'})$. Now $M/H'$ is abelian (as $H/H'\cong C_2^{2n}$) and transitive on each of $U$ and $W$. It follows that the stabiliser $(M/H')_u$ in $M/H'$ of a vertex $u\in U$ fixes $U$ pointwise, and hence is normal in $G^+/H'$. Moreover $(M/H')_u\ne 1$ since $M/H'$ is not semiregular on $V(\Sigma)$, and hence $(M/H')_u$ is transitive on $W$ since $G^+/H'$ is $2$-transitive on $W$. Thus $|M/H'|$ is divisible by $|U|\cdot|W|=2^{2n} = |H/H'|$, and hence $M=H$  so $H/H'=M/H'$ is a minimal normal subgroup of $G/H'$, proving the claim.

\smallskip
Recall that $H/H'$ is regular on  $E(\Sigma_{H'})$, and that  $L:=A(H,X,Y)$ is the stabiliser in $G= H\rtimes L$  of a vertex of $\G$. Since $\G$ is the line graph of $\Sigma$, this implies that $L$ is the stabiliser in $G$ of an edge of $\Sigma$. Moreover, by the Claim, $H/H'$ is a minimal normal subgroup of  $G/H'=(H/H') \rtimes  (LH'/H')$, and hence $(LH')/H'$ is a maximal subgroup of $G/H'$. Thus $G/H'$ is primitive on $E(\Sigma_{H'})$, and hence also its over-group $\Aut(\Sigma)/H'$  is primitive on $E(\Sigma_{H'})$. In these actions, $H/H'\unlhd G/H'$ and  $T/H'=\soc(\Aut(\Sigma)^+/H')\unlhd \Aut(\Sigma)/H'$ are regular normal subgroups, both isomorphic to $C_2^{2n}$. It follows from \cite[Proposition~5.1]{Praeger-PLMS1990} that these regular groups must be equal: thus $H/H'=T/H'$, and so $H=T$ and $H\unlhd\Aut(\Sigma)$. However $\Aut(\Sigma)=\Aut(\G)$, by Lemma~\ref{lem:prop-mixed-dih}(2), and as $G= H\rtimes L$ is the full normaliser of $H$ in $\Aut(\G)$, we conclude that $G=\Aut(\G)$, which contradicts our assumption that $G<\Aut(\G)$. This completes the proof.  \hfill\qed

\section{Proof of Theorem~\ref{non-Cayley}}

The goal of this section is to prove Theorem~\ref{non-Cayley}. Let $X_0=\{x_1,x_2,x_3,x_4\}$, $Y_0=\{y_1,y_2,y_3,y_4\}$, $r$, and $\sigma=(1,2,4,3)$, and let
\[
\P=\lg X_0\cup Y_0\cup\{r\}\mid \mathcal{R}\rg\quad \text{and}\quad \H=\lg X_0\cup Y_0\rg,
\]
where $\mathcal{R}$ is the set of relations given in Definition~\ref{def:H}. Also let $X=\lg X_0\rg$ and $Y=\lg Y_0\rg$. Our first step is to show that $\P$ is a $2$-group and to bound its order.

\begin{lem}\label{lem:1.3a}
\begin{enumerate}
        \item[(a)] $\H\unlhd \P$, $\P=\H\lg r\rg$, and $|r|=8$; and moreover, $x_i^{r^2}=x_{i\sigma}$ and $y_i^{r^2}=y_{i\sigma}$ for $1\leq i\leq 4$.

        \item[(b)] $\H$ has order dividing $2^{152}$, so $\P$ divides $2^{155}$. Moreover
        \begin{enumerate}
            \item[(i)]  $\H/\H'=\lg z\H'\mid z\in X_0\cup Y_0\rg\leq C_2^8$; and
            \item[(ii)] $\H'=\lg [x_i,y_j],[[x_i,y_j],x_k],[[x_i,y_j],y_\ell]\mid 1\leq i,j,k,\ell\leq 4\rg \leq C_2^{144}$.
        \end{enumerate}
\end{enumerate}
\end{lem}

\demo
(a) By Definition~\ref{def:H}, conjugation by $r$ leaves $X_0\cup Y_0$ invariant, and hence $r$ normalises $\H$, so  $\H\unlhd\P$, and $\P=\H\lg r\rg$. Moreover for each $i$ we have
$x_i^{r^2}=y_i^r=x_{i\sigma}$ and $y_i^{r^2}=x_{i\sigma}^r=y_{i\sigma}$. Hence also  $x_i^{r^4}=x_{i\sigma^2}$ and $y_i^{r^4}=y_{i\sigma^2}$ for each $i$, and so since $|\sigma|=4$, we have $r^4\ne 1$. Since $r^8=1$ is a relation in $\R$, it follows that $|r|=8$.

(b) Since $\H=\lg X_0\cup Y_0\rg$, the abelian quotient $\H/\H'=\lg z\H'\mid z\in X_0\cup Y_0\rg$, and since $\R$ contains the relation $z^2=1$ for each of the eight elements $z\in X_0\cup Y_0$, it follows that $\H/\H'\leq C_2^8$, proving part (b)(i). Next, by \cite[Hilfsatz 1.11(a) and (b)]{Huppert},
\[
 \H'=\lg [z,z'],[[z,z'],z'']^h\mid z,z',z''\in X_0\cup Y_0,\ h\in\H\rg.
\]
For all $z,z',z'',z'''\in X_0\cup Y_0$, $\R$ contains the relation $[[[z,z'],z''],z''']=1$, and hence $[[z,z'],z'']\in Z(\H)$.
Also, for all $x,x'\in X_0$ and $y,y'\in Y_0$, $\R$ has relations $[x,x']=[y,y']=[x,y]=1$, and since the inverse $[u,v]^{-1}=[v,u]$ for all $u,v$, it follows that $\H'$ is generated by the set of $16+64+64=144$ commutators given in part (b)(ii). Consider $a=[x_i,y_j], b=[x_{i'},y_{j'}]\in\H'$. Since $\R$ contains the relations $a^2=b^2=x_i^2=y_j^2=1$, it follows that
\[
a^b=[b^{-1}x_ib, b^{-1}y_jb]=[(b^{-1}x_i^{-1}bx_i)x_i, (b^{-1}y_j^{-1}by_j)y_j]=[[b,x_i]x_i,[b,y_j]y_j].
\]
Now $[b,x_i]=[[x_{i'},y_{j'}],x_i]\in Z(\H)$ and $[b,y_j]=[[x_{i'},y_{j'}],y_j]\in Z(\H)$, and hence the last expression for $a^b$ is equal to $[x_i,y_j]=a$. Thus $a^b=a$, that is, $a$ and $b$ commute. It follows that $\H'$ is abelian. Finally, since $\R$ contains the relations $[x_i,y_j]^2=[[x_i,y_j],x_k]^2=[[x_i,y_j],y_\ell]^2=1$, we conclude that $\H'$ is an elementary abelian $2$-group generated by the given $144$ generators in (b)(ii). Hence $\H'\leq C_2^{144}$, so we have shown that $\P$ is a $2$-group, and $|\P|\leq 8\cdot|\H|\leq 8\cdot|\H/\H'|\cdot|\H'|\leq 2^{3+8+144}$. This proves part (b).
\qed

\subsection{A Magma computation}\label{sub:magma}
From now on we make use of a Magma~\cite{BCP} computation. We have made available the Magma programs in an Appendix in Section~\ref{s:magma} of this paper.

First, the group $\P$ is input in the category GrpFP via the presentation given in Definition~\ref{def:H}. Since by Lemma~\ref{lem:1.3a}, $|\P|\leq 2^{155}$, $155$ is an upper bound for the length of the lower exponent $2$ central series of $\P$. We used the {\tt pQuotient} command to construct the largest $2$-quotient $P$ of $\P$ having lower exponent-$2$ class at most 200 as a group in the category GrpPC. We are guaranteed to have $\P\cong P$. Our computation in Magma~\cite{BCP} showed that $|\P|=8\cdot|\H|=2^{59}$. Hence $\P=\H\rtimes\lg r\rg$.

We next used Magma~\cite{BCP} to compute the derived quotient of $\H$ and found that $\H/\H'\cong C_2^{8}$, and that $X=\lg X_0\rg\cong C_2^4$ and $Y=\lg Y_0\rg\cong C_2^4$. Thus:

\begin{lem}\label{lem:1.30}
$\H$ is a $4$-dimensional mixed dihedral group of order $2^{56}$ relative to $X$ and $Y$, and Theorem~\ref{non-Cayley}~(1) holds.
\end{lem}

\smallskip
We now define several maps on $X_0\cup Y_0$ as follows:
\[
\begin{array}{l}
\a_1: x_1\mapsto x_1x_2, x_2\mapsto x_2x_3, x_3\mapsto x_3x_4, x_4\mapsto x_1x_2x_3, y_i\mapsto y_i (i=1,2,3,4),\\

\a_2: x_i\mapsto x_i (i=1,2,3,4), y_1\mapsto y_1y_2, y_2\mapsto y_2y_3, y_3\mapsto y_3y_4, y_4\mapsto y_1y_2y_3,\\

\b_1: x_1\mapsto x_2, x_2\mapsto x_3, x_3\mapsto x_4, x_4\mapsto x_1x_2x_3x_4, y_i\mapsto y_i (i=1,2,3,4),\\

\b_2: x_i\mapsto x_i (i=1,2,3,4), y_1\mapsto y_2, y_2\mapsto y_3, y_3\mapsto y_4, y_4\mapsto y_1y_2y_3y_4. \\


\end{array}
\]
Using the {\tt hom} command in Magma~\cite{BCP}, we verified that each of these four maps induces an automorphism of $\H$ which leaves invariant each of $X$ and $Y$. In the following, we also use $\a_1,\a_2,\b_1$ and $\b_2$ to denote the corresponding automorphisms of $\H$ induced by the above four maps, respectively.

It is easy to check that $\beta_i$ has order $5$ for $i=1, 2$. Next, using the facts that $X\cong Y\cong C_2^4$ and the definitions of these maps, we check (by hand) that
\[
\begin{array}{l}
\a_1^3: x_1\mapsto x_1x_2x_3x_4, x_2\mapsto x_1, x_3\mapsto x_2, x_4\mapsto x_3, y_i\mapsto y_i (i=1,2,3,4),\\
\a_2^3: x_i\mapsto x_i (i=1,2,3,4), y_1\mapsto y_1y_2y_3y_4, y_2\mapsto y_1, y_3\mapsto y_2, y_4\mapsto y_3.
\end{array}\]
This implies that, for each $i=1$ or $2$,  $\a_i^3=\b_i^{-1}$ and hence $\alpha_i^{15}$ is the identity. We have seen that $\a_i^3\ne 1$ and it is easy to see that $\a_i^5\ne 1$ (for example, $\a_1^5: x_3\mapsto x_2x_4$), and hence $\a_1$ and $\a_2$ have order $15$. Moreover, it is clear that $\a_1$ and $\a_2$ commute in their actions on $\H$.

As we noted above, and in Lemma~\ref{lem:1.3a}, $|r|=8$, and the conjugation action of $r$ on $\H$ is an automorphism of $\H$ of order $8$. Its action on the generators is given by relations in $\R$, namely
\[
r: x_1\mapsto y_1, x_2\mapsto y_2, x_3\mapsto y_3, x_4\mapsto y_4, y_1\mapsto x_2, y_2\mapsto x_4, y_3\mapsto x_1, y_4\mapsto x_3
\]
and the action of $r^2$ is recorded
in Lemma~\ref{lem:1.3a}(a). In particular $r$ interchanges $X$ and $Y$.
Further we may check that the three automorphisms $r, \a_1$ and $\a_2$ of $\H$ satisfy relations: $\a_1^r=\a_2$ and $\a_2^r=\a_1^2$.
We summarise this information in the next lemma.

\begin{lem}\label{lem:1.3b}
Let $K:=\lg \a_1, \a_2, r\rg \leq \Aut(\H)$, and let $A=\H\rtimes K$ with the actions above. Then $K=(\lg\a_1\rg\times\lg\a_2\rg)\rtimes \lg r\rg\cong (C_{15}\times C_{15})\rtimes C_8$, $K$ leaves $S =(X\cup Y)\setminus\{1\}$ invariant, its subgroup $\lg \a_1, \a_2, r^2\rg$ leaves each of $X$ and $Y$ invariant while  $r$ interchanges $X$ and $Y$, and $\P$ is a Sylow $2$-subgroup of $A$.
\end{lem}


\subsection{The graphs $\G=\G(\H, X, Y)$ and $\Sigma=\Sigma(\H, X, Y)$}

As we observed in Subsection~\ref{sub:magma}, $\H$ is a mixed dihedral group relative to $X$ and $Y$, and we let $\G=\G(\H, X, Y)$ and $\Sigma=\Sigma(\H, X, Y)$ be the graphs as in Subsection~\ref{s:mix-dih}.
Thus $\G$ is the Cayley graph $\Cay(\H, S)$, where
\[
\begin{array}{l}
S=(X\cup Y)\setminus\{1\},
\end{array}
\]
and $\Sigma$ has vertex set and edge set as follows:
\[\begin{array}{l}
V=\{Xh, Yg: h,g\in \H\},\\
E=\{\{Xh, Yg\}:\   h,g\in H\ \mbox{and}\ Xh\cap Yg\neq\emptyset\}. 
\end{array}
\]
By Lemma~\ref{lem:prop-mixed-dih}~(1) and (2), $\Gamma$ is isomorphic to the line graph of $\Sigma$ and  $\Aut(\G)=\Aut(\Sigma)$.
By Lemma~\ref{lem:1.3b}, $K=\lg \a_1,\a_2,r\rg$ fixes $S$ setwise, and hence $K$ is contained in the group $A(\H, X, Y)$ in \eqref{eq-ahxy}. We prove that equality holds, again with some help from Magma~\cite{BCP}.

\begin{lem}\label{lem:1.3c}
The group $K=\lg \a_1, \a_2, r\rg \leq \Aut(\H)$ is equal to $A(\H,X,Y)$. Moreover $K=(C_{15}\times C_{15})\rtimes C_8$, $K$ is transitive on $S$, and $K_{(X)}\ne1$.
\end{lem}

\demo
By Lemma~\ref{lem:prop-mixed-dih}~(4), $A(\H, X, Y)\leq(\Aut(X)\times\Aut(Y))\rtimes C_2 \cong(\GL(4,2)\times\GL(4,2))\rtimes C_2,$ where the $C_2$ interchanges $X$ and $Y$. By Lemma~\ref{lem:1.3b}, $r\in K\leq A(\H,X,Y)$ interchanges $X$ and $Y$, so the subgroup $L$ of $A(\H,X,Y)$ fixing each of $X, Y$ setwise has index $2$, and $K\cap L = \lg \a_1, \a_2, r^2\rg$. For $\Delta\in\{X, Y, S\}$ and any $g\in L$ and $B\leq L$, let $g^\Delta$ and $B^{\Delta}$ denote the induced action on $\Delta$. Then
\begin{equation}\label{eq-KLX}
    (K\cap L)^{X}=\lg(\a_1)^{X}\rg\rtimes \lg(r^2)^{X}\rg\cong C_{15}\rtimes C_4,
\end{equation}
which is an irreducible subgroup of $\Aut(X)\cong\GL(4,2)$ acting transitively on $X\setminus\{1\}$, so $K$ is transitive on $S$.
Also (since by Lemma~\ref{lem:1.3a}, $\lg r^2\rg$ acts faithfully on $X$), the pointwise stabiliser $(K\cap L)_{(X)}=\lg \a_2\rg\cong C_{15}$, and acts transitively on $Y\setminus\{1\}$. Since $S$ generates $\H$, the group $A(\H,X,Y)$ acts faithfully on $S$, and we have $C_{15}\times C_{15} = \lg\a_1\rg^X\times\lg\a_2\rg^Y \leq L^X\times L^Y$, with $r$ interchanging $L^X$ and $L^Y$. In particular $L^Y\cong L^X$. We note that the only group $B$ such that  $(K\cap L)^{X} < B<\GL(4,2)$  is $(A_5\times C_3)\rtimes C_2$ (see \cite[p.22]{Atlas}).

Suppose first that $L$ is insoluble. Then  $L^X (\cong L^Y) = (A_5\times C_3)\rtimes C_2$ or $\GL(4,2)$. Since $(L_{(X)})^Y$ is a normal subgroup of $L^{Y}$ containing $\a_2^Y$, it follows that $L_{(X)}\cong C_3\times A_5, (A_5\times C_3)\rtimes C_2$, or $\GL(4,2)$, and moreover that  either $L\cong \Aut(X)\times \Aut(Y)$, or
$\soc(L)=\soc(L_{(X)})\times\soc(L_{(Y)})=(A_5\times C_3)^2$. In the latter case $\soc(L)$ has centre
$Z:=\lg\a_1^5\rg\times\lg \a_2^5\rg\cong C_3\times C_3$ and $\soc(L)$ is equal to the full centraliser in $\Aut(X)\times \Aut(Y)$ of $Z$. Thus in either case $L$ contains the full centraliser of $Z$.
It is easy to check that
\[
\begin{array}{l}
\a_1^5: x_1\mapsto x_1x_3x_4, x_2\mapsto x_1x_3, x_3\mapsto x_2x_4, x_4\mapsto x_1x_2x_4, y_i\mapsto y_i (i=1,2,3,4),\\
\a_2^5: x_i\mapsto x_i (i=1,2,3,4), y_1\mapsto y_1y_3y_4, y_2\mapsto y_1y_3, y_3\mapsto y_2y_4, y_4\mapsto y_1y_2y_4.
\end{array}
\]
Consider the element $\s_1\in \Aut(X)\times\Aut(Y)$ defined by
\[\begin{array}{l}
\s_1: x_1\mapsto x_1, x_2\mapsto x_1x_3, x_3\mapsto x_2x_3, x_4\mapsto x_2x_4, y_i\mapsto y_i \ (i=1,2,3,4), \\

\end{array}
\]
Now $\a_1^5$ and $\s_1$ both fix $Y$ pointwise and the following equalities can easily be verified:
\[
\begin{array}{l}
{x_1^{\s_1\a_1^5}=x_1x_3x_4=(x_1)^{\a_1^5\s_1}, x_2^{\s_1\a_1^5}=x_1x_2x_3=x_2^{\a_1^5\s_1},} \\

{x_3^{\s_1\a_1^5}=x_1x_2x_3x_4=x_3^{\a_1^5\s_1}, x_4^{\s_1\a_1^5}=x_2x_3x_4=x_4^{\a_1^5\s_1}.}\\
\end{array}
\]
This implies that $\s_1\a_1^5=\a_1^5\s_1$. Also $\s_1\a_2^5=\a_2^5\s_1$ since $\s_1, \a_2^5$ act trivially on $Y, X$ respectively. Thus $\s_1$ centralises $Z$ and hence lies in $L$. Since $L\leq A(\H, X, Y)$ acts faithfully on $S$, it follows $\s_1$ can be extended to an element in $A(\H, X, Y)$.
However, use of the {\tt hom} command in Magma~\cite{BCP} shows that $\s_1$ does not induce an automorphism of $\H$, which is a contradiction.

Thus, $L$ is soluble, so both $L^{X}$ and $L^{Y}$ are soluble. Since as we observed above the only proper overgroups of $(K\cap L)^X$ in $\GL(4,2)$ are insoluble, it follows that $L^X=(K\cap L)^X$ as in \eqref{eq-KLX} and similarly $(K\cap L)^{Y} = \lg(\a_2)^{Y}\rg\rtimes \lg(r^2)^{Y}\rg \cong C_{15}\rtimes C_4$. Thus
\[
L\cong L^S\leq L^{X}\times L^{Y} = (\lg(\a_1)^{X}\rg\rtimes \lg(r^2)^{X}\rg) \times (\lg(\a_2)^{Y}\rg\rtimes \lg(r^2)^{Y}\rg)
\]
and $L^S$ contains $(K\cap L)^S = (\lg \a_1^{X}\rg\times\lg \a_2^{Y}\rg)\rtimes \lg (r^2)^S\rg$.
Therefore, if $L$ is strictly larger than $K\cap L$, then $L\cap (\lg(r^2)^{X}\rg \times \lg(r^2)^{Y}\rg)$ must contain $\lg(r^4)^{X}\rg \times \lg(r^4)^{Y}\rg$. In particular, in this case, $L$ must contain $(r^4)^{X}$. Now $(r^4)^{X}$ fixes $Y$ pointwise and, by Lemma~\ref{lem:1.3a}, maps $x_i\mapsto x_{i\s^2}$ for each $i$, with $\s^2=(1,4)(2,3)$. However, use of the {\tt hom} command in Magma~\cite{BCP} shows that $(r^4)^{X}$ does not induce an automorphism of $\H$, and we have a contradiction. Thus we conclude that $L=K\cap L$, and since $|K:K\cap L|=|A(\H,X,Y):L|=2$, this implies that  $A(\H, X, Y)=K$, completing the proof.
\hfill\qed

Next we determine the full automorphism group of $\G$.

\begin{lem}\label{lem:1.3d}
$\Aut(\Sigma)=\Aut(\G)=\H\rtimes A(\H,X,Y)$. Moreover $\G$ is normal Cayley graph and is edge-transitive, and Theorem~\ref{non-Cayley}~(2) holds.
\end{lem}

\demo
By Lemma~\ref{lem:1.3c}, the group $A(\H, X, Y)=K$ is transitive on $S$, and hence $\H\rtimes K$ is edge-transitive on $\G$. Moreover, the pointwise stabiliser $K_{(X)}$ is nontrivial and $K$ does not contain $\Aut(X)\times\Aut(Y)$. Hence, by Lemma~\ref{lem:1.30} the conditions of Theorem~\ref{mix-dihedrant-characterization} all hold, and we conclude  that $\Aut(\G)=\H\rtimes A(\H, X, Y)$. Thus $\G$ is a normal Cayley graph. Thus Theorem~\ref{non-Cayley}~(2) holds since $\Aut(\Sigma)=\Aut(\G)$.
\hfill \qed

Now we complete the proof of Theorem~\ref{non-Cayley}. Parts (1) and (2) were proved in Lemmas~\ref{lem:1.30} and~\ref{lem:1.3d}. We now prove part (3). By Lemma~\ref{lem:1.30},  $|\H|=2^{56}$, and as $\G$ is a Cayley graph for $\H$, $|V(\G)|=2^{56}$. Also, by Lemma~\ref{lem:prop-mixed-dih}~(3), $\H$ has two orbits in $V(\Sigma)$, with vertex stabilisers $X, Y$ of the vertices $X, Y$ respectively. Hence $|V(\Sigma)|=2\cdot |\H:X|= 2^{1+56-4}=2^{53}$.  By Lemma~\ref{lem:1.3d}, $\G$ is $\Aut(\G)$-edge-transitive, and by Lemma~\ref{lem:prop-mixed-dih}~(1), $\G$ is isomorphic to the line graph of $\Sigma$. Thus all the conditions of Theorem~\ref{iff-characterization}~(b) hold (with the group $G=\Aut(\G)$), and it follows that  $\Sigma$ is $(\Aut(\G),2)$-arc-transitive, and is an $\H'$-normal cover, of order $2^{53}$, of $\K_{2^4,2^4}$.


It remains to prove that $\Sigma$ is not a Cayley graph. Suppose to the contrary that $\Sigma$ is a Cayley graph.
Then $\Aut(\G)=\H: K$ has a subgroup, say $G$, acting regularly on $V(\Sigma)$. In particular $|G|=2^{53}$. By Lemma~\ref{lem:1.3b}, $\P=\H\rtimes \lg r\rg$ is a Sylow $2$-subgroup of $G$, and $|\P|=|\H|\cdot 2^3$. By Sylow's theorem, we may assume that $G\leq \P$. Note that $X$ is a vertex of $\Sigma$ and the (setwise) stabiliser of $X$ in $\P$ is
\[
\P_X=(\H\rtimes\lg r\rg)_X=X\rtimes\lg r^2\rg.
\]
As $G$ is regular on $V(\Sigma)$,  $|G\cap (\H\rtimes\lg r\rg)_{X}|=1$. However, by using Magma~\cite{BCP}, we see that $\P=\H\rtimes\lg r\rg$ has no subgroup $T$ of order $2^{53}$ such that $|T\cap (\H\rtimes \lg r\rg)_{X}|=1$, and we have a contradiction. Thus $\Sigma$ is not a Cayley graph, and the proof of part (3), and hence of Theorem~\ref{non-Cayley}, is complete.\hfill\qed

\section{Appendix: Magma programs for proving Theorem~\ref{non-Cayley}}\label{s:magma}

\f{Input the group $\P$:}\smallskip
\\
{\tt G<x1,x2,x3,x4,y1,y2,y3,y4,r>:=Group<x1,x2,x3,x4,y1,y2,y3,y4,r|
\\
x1\textasciicircum2, x2\textasciicircum2, x3\textasciicircum2, x4\textasciicircum2,
y1\textasciicircum2, y2\textasciicircum2, y3\textasciicircum2, y4\textasciicircum2, r\textasciicircum8,
\\
x1\textasciicircum r=y1,x2\textasciicircum r=y2,x3\textasciicircum r=y3,x4\textasciicircum r=y4,
y1\textasciicircum r=x2,y2\textasciicircum r=x4,y3\textasciicircum r=x1,y4\textasciicircum r=x3,
\\
(x1,x2)=(x1,x3)=(x1,x4)=(x2,x3)=(x2,x4)=(x3,x4)=1,
\\
(y1,y2)=(y1,y3)=(y1,y4)=(y2,y3)=(y2,y4)=(y3,y4)=1,
\\
(x1,y1)\textasciicircum2=1,(x1,y2)\textasciicircum2=1,(x1,y3)\textasciicircum2=1,(x1,y4)\textasciicircum2=1,
\\
(x2,y1)\textasciicircum2=1,(x2,y2)\textasciicircum2=1,(x2,y3)\textasciicircum2=1,(x2,y4)\textasciicircum2=1,
\\
(x3,y1)\textasciicircum2=1,(x3,y2)\textasciicircum2=1,(x3,y3)\textasciicircum2=1,(x3,y4)\textasciicircum2=1,
\\
(x4,y1)\textasciicircum2=1,(x4,y2)\textasciicircum2=1,(x4,y3)\textasciicircum2=1,(x4,y4)\textasciicircum2=1,
\\
((x1,y1),x2)\textasciicircum2=((x1,y1),x3)\textasciicircum2=((x1,y1),x4)\textasciicircum2=((x1,y1),y2)\textasciicircum2=
\\
((x1,y1),y3)\textasciicircum2=((x1,y1),y4)\textasciicircum2=1,
\\
((x1,y2),x2)\textasciicircum2=((x1,y2),x3)\textasciicircum2=((x1,y2),x4)\textasciicircum2=((x1,y2),y1)\textasciicircum2=
\\
((x1,y2),y3)\textasciicircum2=((x1,y2),y4)\textasciicircum2=1,
\\
((x1,y3),x2)\textasciicircum2=((x1,y3),x3)\textasciicircum2=((x1,y3),x4)\textasciicircum2=((x1,y3),y2)\textasciicircum2=
\\
((x1,y3),y1)\textasciicircum2=((x1,y3),y4)\textasciicircum2=1,
\\
((x1,y4),x2)\textasciicircum2=((x1,y4),x3)\textasciicircum2=((x1,y4),x4)\textasciicircum2=((x1,y4),y2)\textasciicircum2=
\\
((x1,y4),y3)\textasciicircum2=((x1,y4),y1)\textasciicircum2=1,
\\
((x2,y1),x1)\textasciicircum2=((x2,y1),x3)\textasciicircum2=((x2,y1),x4)\textasciicircum2=((x2,y1),y2)\textasciicircum2=
\\
((x2,y1),y3)\textasciicircum2=((x2,y1),y4)\textasciicircum2=1,
\\
((x2,y2),x1)\textasciicircum2=((x2,y2),x3)\textasciicircum2=((x2,y2),x4)\textasciicircum2=((x2,y2),y1)\textasciicircum2=
\\
((x2,y2),y3)\textasciicircum2=((x2,y2),y4)\textasciicircum2=1,
\\
((x2,y3),x1)\textasciicircum2=((x2,y3),x3)\textasciicircum2=((x2,y3),x4)\textasciicircum2=((x2,y3),y2)\textasciicircum2=
\\
((x2,y3),y1)\textasciicircum2=((x2,y3),y4)\textasciicircum2=1,
\\
((x2,y4),x1)\textasciicircum2=((x2,y4),x3)\textasciicircum2=((x2,y4),x4)\textasciicircum2=((x2,y4),y2)\textasciicircum2=
\\
((x2,y4),y3)\textasciicircum2=((x2,y4),y1)\textasciicircum2=1,
\\
((x3,y1),x2)\textasciicircum2=((x3,y1),x1)\textasciicircum2=((x3,y1),x4)\textasciicircum2=((x3,y1),y2)\textasciicircum2=
\\
((x3,y1),y3)\textasciicircum2=((x3,y1),y4)\textasciicircum2=1,
\\
((x3,y2),x2)\textasciicircum2=((x3,y2),x1)\textasciicircum2=((x3,y2),x4)\textasciicircum2=((x3,y2),y1)\textasciicircum2=
\\
((x3,y2),y3)\textasciicircum2=((x3,y2),y4)\textasciicircum2=1,
\\
((x3,y3),x2)\textasciicircum2=((x3,y3),x1)\textasciicircum2=((x3,y3),x4)\textasciicircum2=((x3,y3),y2)\textasciicircum2=
\\
((x3,y3),y1)\textasciicircum2=((x3,y3),y4)\textasciicircum2=1,
\\
((x3,y4),x2)\textasciicircum2=((x3,y4),x1)\textasciicircum2=((x3,y4),x4)\textasciicircum2=((x3,y4),y2)\textasciicircum2=
\\
((x3,y4),y3)\textasciicircum2=((x3,y4),y1)\textasciicircum2=1,
\\
((x4,y1),x2)\textasciicircum2=((x4,y1),x1)\textasciicircum2=((x4,y1),x3)\textasciicircum2=((x4,y1),y2)\textasciicircum2=
\\
((x4,y1),y3)\textasciicircum2=((x4,y1),y4)\textasciicircum2=1,
\\
((x4,y2),x2)\textasciicircum2=((x4,y2),x1)\textasciicircum2=((x4,y2),x3)\textasciicircum2=((x4,y2),y1)\textasciicircum2=
\\
((x4,y2),y3)\textasciicircum2=((x4,y2),y4)\textasciicircum2=1,
\\
((x4,y3),x2)\textasciicircum2=((x4,y3),x1)\textasciicircum2=((x4,y3),x3)\textasciicircum2=((x4,y3),y2)\textasciicircum2=
\\
((x4,y3),y1)\textasciicircum2=((x4,y3),y4)\textasciicircum2=1,
\\
((x4,y4),x2)\textasciicircum2=((x4,y4),x1)\textasciicircum2=((x4,y4),x3)\textasciicircum2=((x4,y4),y2)\textasciicircum2=
\\
((x4,y4),y3)\textasciicircum2=((x4,y4),y1)\textasciicircum2=1,
\\
(x1,((x1,y1),x2))=(x2,((x1,y1),x2))=(x3,((x1,y1),x2))=(x4,((x1,y1),x2))=
\\
(y1,((x1,y1),x2))=(y2,((x1,y1),x2))=(y3,((x1,y1),x2))=(y4,((x1,y1),x2))=1,
\\
(x1,((x1,y1),x3))=(x2,((x1,y1),x3))=(x3,((x1,y1),x3))=(x4,((x1,y1),x3))=
\\
(y1,((x1,y1),x3))=(y2,((x1,y1),x3))=(y3,((x1,y1),x3))=(y4,((x1,y1),x3))=1,
\\
(x1,((x1,y1),x4))=(x2,((x1,y1),x4))=(x3,((x1,y1),x4))=(x4,((x1,y1),x4))=
\\
(y1,((x1,y1),x4))=(y2,((x1,y1),x4))=(y3,((x1,y1),x4))=(y4,((x1,y1),x4))=1,
\\
(x1,((x1,y1),y2))=(x2,((x1,y1),y2))=(x3,((x1,y1),y2))=(x4,((x1,y1),y2))=
\\
(y1,((x1,y1),y2))=(y2,((x1,y1),y2))=(y3,((x1,y1),y2))=(y4,((x1,y1),y2))=1,
\\
(x1,((x1,y1),y3))=(x2,((x1,y1),y3))=(x3,((x1,y1),y3))=(x4,((x1,y1),y3))=
\\
(y1,((x1,y1),y3))=(y2,((x1,y1),y3))=(y3,((x1,y1),y3))=(y4,((x1,y1),y3))=1,
\\
(x1,((x1,y1),y4))=(x2,((x1,y1),y4))=(x3,((x1,y1),y4))=(x4,((x1,y1),y4))=
\\
(y1,((x1,y1),y4))=(y2,((x1,y1),y4))=(y3,((x1,y1),y4))=(y4,((x1,y1),y4))=1,
\\
(x1,((x1,y2),x2))=(x2,((x1,y2),x2))=(x3,((x1,y2),x2))=(x4,((x1,y2),x2))=
\\
(y1,((x1,y2),x2))=(y2,((x1,y2),x2))=(y3,((x1,y2),x2))=(y4,((x1,y2),x2))=1,
\\
(x1,((x1,y2),x3))=(x2,((x1,y2),x3))=(x3,((x1,y2),x3))=(x4,((x1,y2),x3))=
\\
(y1,((x1,y2),x3))=(y2,((x1,y2),x3))=(y3,((x1,y2),x3))=(y4,((x1,y2),x3))=1,
\\
(x1,((x1,y2),x4))=(x2,((x1,y2),x4))=(x3,((x1,y2),x4))=(x4,((x1,y2),x4))=
\\
(y1,((x1,y2),x4))=(y2,((x1,y2),x4))=(y3,((x1,y2),x4))=(y4,((x1,y2),x4))=1,
\\
(x1,((x1,y2),y1))=(x2,((x1,y2),y1))=(x3,((x1,y2),y1))=(x4,((x1,y2),y1))=
\\
(y1,((x1,y2),y1))=(y2,((x1,y2),y1))=(y3,((x1,y2),y1))=(y4,((x1,y2),y1))=1,
\\
(x1,((x1,y2),y3))=(x2,((x1,y2),y3))=(x3,((x1,y2),y3))=(x4,((x1,y2),y3))=
\\
(y1,((x1,y2),y3))=(y2,((x1,y2),y3))=(y3,((x1,y2),y3))=(y4,((x1,y2),y3))=1,
\\
(x1,((x1,y2),y4))=(x2,((x1,y2),y4))=(x3,((x1,y2),y4))=(x4,((x1,y2),y4))=
\\
(y1,((x1,y2),y4))=(y2,((x1,y2),y4))=(y3,((x1,y2),y4))=(y4,((x1,y2),y4))=1,
\\
(x1,((x1,y3),x2))=(x2,((x1,y3),x2))=(x3,((x1,y3),x2))=(x4,((x1,y3),x2))=
\\
(y1,((x1,y3),x2))=(y2,((x1,y3),x2))=(y3,((x1,y3),x2))=(y4,((x1,y3),x2))=1,
\\
(x1,((x1,y3),x3))=(x2,((x1,y3),x3))=(x3,((x1,y3),x3))=(x4,((x1,y3),x3))=
\\
(y1,((x1,y3),x3))=(y2,((x1,y3),x3))=(y3,((x1,y3),x3))=(y4,((x1,y3),x3))=1,
\\
(x1,((x1,y3),x4))=(x2,((x1,y3),x4))=(x3,((x1,y3),x4))=(x4,((x1,y3),x4))=
\\
(y1,((x1,y3),x4))=(y2,((x1,y3),x4))=(y3,((x1,y3),x4))=(y4,((x1,y3),x4))=1,
\\
(x1,((x1,y3),y2))=(x2,((x1,y3),y2))=(x3,((x1,y3),y2))=(x4,((x1,y3),y2))=
\\
(y1,((x1,y3),y2))=(y2,((x1,y3),y2))=(y3,((x1,y3),y2))=(y4,((x1,y3),y2))=1,
\\
(x1,((x1,y3),y1))=(x2,((x1,y3),y1))=(x3,((x1,y3),y1))=(x4,((x1,y3),y1))=
\\
(y1,((x1,y3),y3))=(y2,((x1,y3),y1))=(y3,((x1,y3),y1))=(y4,((x1,y3),y1))=1,
\\
(x1,((x1,y3),y4))=(x2,((x1,y3),y4))=(x3,((x1,y3),y4))=(x4,((x1,y3),y4))=
\\
(y1,((x1,y3),y4))=(y2,((x1,y3),y4))=(y3,((x1,y3),y4))=(y4,((x1,y3),y4))=1,
\\
(x1,((x1,y4),x2))=(x2,((x1,y4),x2))=(x3,((x1,y4),x2))=(x4,((x1,y4),x2))=
\\
(y1,((x1,y4),x2))=(y2,((x1,y4),x2))=(y3,((x1,y4),x2))=(y4,((x1,y4),x2))=1,
\\
(x1,((x1,y4),x3))=(x2,((x1,y4),x3))=(x3,((x1,y4),x3))=(x4,((x1,y4),x3))=
\\
(y1,((x1,y4),x3))=(y2,((x1,y4),x3))=(y3,((x1,y4),x3))=(y4,((x1,y4),x3))=1,
\\
(x1,((x1,y4),x4))=(x2,((x1,y4),x4))=(x3,((x1,y4),x4))=(x4,((x1,y4),x4))=
\\
(y1,((x1,y4),x4))=(y2,((x1,y4),x4))=(y3,((x1,y4),x4))=(y4,((x1,y4),x4))=1,
\\
(x1,((x1,y4),y2))=(x2,((x1,y4),y2))=(x3,((x1,y4),y2))=(x4,((x1,y4),y2))=
\\
(y1,((x1,y4),y2))=(y2,((x1,y4),y2))=(y3,((x1,y4),y2))=(y4,((x1,y4),y2))=1,
\\
(x1,((x1,y4),y3))=(x2,((x1,y4),y3))=(x3,((x1,y4),y3))=(x4,((x1,y4),y3))=
\\
(y1,((x1,y4),y3))=(y2,((x1,y4),y3))=(y3,((x1,y4),y3))=(y4,((x1,y4),y3))=1,
\\
(x1,((x1,y4),y1))=(x2,((x1,y4),y1))=(x3,((x1,y4),y1))=(x4,((x1,y4),y1))=
\\
(y1,((x1,y4),y1))=(y2,((x1,y4),y1))=(y3,((x1,y4),y1))=(y4,((x1,y4),y1))=1,
\\
(x1,((x2,y1),x1))=(x2,((x2,y1),x1))=(x3,((x2,y1),x1))=(x4,((x2,y1),x1))=
\\
(y1,((x2,y1),x1))=(y2,((x2,y1),x1))=(y3,((x2,y1),x1))=(y4,((x2,y1),x1))=1,
\\
(x1,((x2,y1),x3))=(x2,((x2,y1),x3))=(x3,((x2,y1),x3))=(x4,((x2,y1),x3))=
\\
(y1,((x2,y1),x3))=(y2,((x2,y1),x3))=(y3,((x2,y1),x3))=(y4,((x2,y1),x3))=1,
\\
(x1,((x2,y1),x4))=(x2,((x2,y1),x4))=(x3,((x2,y1),x4))=(x4,((x2,y1),x4))=
\\
(y1,((x2,y1),x4))=(y2,((x2,y1),x4))=(y3,((x2,y1),x4))=(y4,((x2,y1),x4))=1,
\\
(x1,((x2,y1),y2))=(x2,((x2,y1),y2))=(x3,((x2,y1),y2))=(x4,((x2,y1),y2))=
\\
(y1,((x2,y1),y2))=(y2,((x2,y1),y2))=(y3,((x2,y1),y2))=(y4,((x2,y1),y2))=1,
\\
(x1,((x2,y1),y3))=(x2,((x2,y1),y3))=(x3,((x2,y1),y3))=(x4,((x2,y1),y3))=
\\
(y1,((x2,y1),y3))=(y2,((x2,y1),y3))=(y3,((x2,y1),y3))=(y4,((x2,y1),y3))=1,
\\
(x1,((x2,y1),y4))=(x2,((x2,y1),y4))=(x3,((x2,y1),y4))=(x4,((x2,y1),y4))=
\\
(y1,((x2,y1),y4))=(y2,((x2,y1),y4))=(y3,((x2,y1),y4))=(y4,((x2,y1),y4))=1,
\\
(x1,((x2,y2),x1))=(x2,((x2,y2),x1))=(x3,((x2,y2),x1))=(x4,((x2,y2),x1))=
\\
(y1,((x2,y2),x1))=(y2,((x2,y2),x1))=(y3,((x2,y2),x1))=(y4,((x2,y2),x1))=1,
\\
(x1,((x2,y2),x3))=(x2,((x2,y2),x3))=(x3,((x2,y2),x3))=(x4,((x2,y2),x3))=
\\
(y1,((x2,y2),x3))=(y2,((x2,y2),x3))=(y3,((x2,y2),x3))=(y4,((x2,y2),x3))=1,
\\
(x1,((x2,y2),x4))=(x2,((x2,y2),x4))=(x3,((x2,y2),x4))=(x4,((x2,y2),x4))=
\\
(y1,((x2,y2),x4))=(y2,((x2,y2),x4))=(y3,((x2,y2),x4))=(y4,((x2,y2),x4))=1,
\\
(x1,((x2,y2),y1))=(x2,((x2,y2),y1))=(x3,((x2,y2),y1))=(x4,((x2,y2),y1))=
\\
(y1,((x2,y2),y1))=(y2,((x2,y2),y1))=(y3,((x2,y2),y1))=(y4,((x2,y2),y1))=1,
\\
(x1,((x2,y2),y3))=(x2,((x2,y2),y3))=(x3,((x2,y2),y3))=(x4,((x2,y2),y3))=
\\
(y1,((x2,y2),y3))=(y2,((x2,y2),y3))=(y3,((x2,y2),y3))=(y4,((x2,y2),y3))=1,
\\
(x1,((x2,y2),y4))=(x2,((x2,y2),y4))=(x3,((x2,y2),y4))=(x4,((x2,y2),y4))=
\\
(y1,((x2,y2),y4))=(y2,((x2,y2),y4))=(y3,((x2,y2),y4))=(y4,((x2,y2),y4))=1,
\\
(x1,((x2,y3),x1))=(x2,((x2,y3),x1))=(x3,((x2,y3),x1))=(x4,((x2,y3),x1))=
\\
(y1,((x2,y3),x1))=(y2,((x2,y3),x1))=
(y3,((x2,y3),x1))=(y4,((x2,y3),x1))=1,
\\
(x1,((x2,y3),x3))=(x2,((x2,y3),x3))=(x3,((x2,y3),x3))=(x4,((x2,y3),x3))=
\\
(y1,((x2,y3),x3))=(y2,((x2,y3),x3))=(y3,((x2,y3),x3))=(y4,((x2,y3),x3))=1,
\\
(x1,((x2,y3),x4))=(x2,((x2,y3),x4))=(x3,((x2,y3),x4))=(x4,((x2,y3),x4))=
\\
(y1,((x2,y3),x4))=(y2,((x2,y3),x4))=(y3,((x2,y3),x4))=(y4,((x2,y3),x4))=1,
\\
(x1,((x2,y3),y2))=(x2,((x2,y3),y2))=(x3,((x2,y3),y2))=(x4,((x2,y3),y2))=
\\
(y1,((x2,y3),y2))=(y2,((x2,y3),y2))=(y3,((x2,y3),y2))=(y4,((x2,y3),y2))=1,
\\
(x1,((x2,y3),y1))=(x2,((x2,y3),y1))=(x3,((x2,y3),y1))=(x4,((x2,y3),y1))=
\\
(y1,((x2,y3),y1))=(y2,((x2,y3),y1))=(y3,((x2,y3),y1))=(y4,((x2,y3),y1))=1,
\\
(x1,((x2,y3),y4))=(x2,((x2,y3),y4))=(x3,((x2,y3),y4))=(x4,((x2,y3),y4))=
\\
(y1,((x2,y3),y4))=(y2,((x2,y3),y4))=(y3,((x2,y3),y4))=(y4,((x2,y3),y4))=1,
\\
(x1,((x2,y4),x1))=(x2,((x2,y4),x1))=(x3,((x2,y4),x1))=(x4,((x2,y4),x1))=
\\
(y1,((x2,y4),x1))=(y2,((x2,y4),x1))=(y3,((x2,y4),x1))=(y4,((x2,y4),x1))=1,
\\
(x1,((x2,y4),x3))=(x2,((x2,y4),x3))=(x3,((x2,y4),x3))=(x4,((x2,y4),x3))=
\\
(y1,((x2,y4),x3))=(y2,((x2,y4),x3))=(y3,((x2,y4),x3))=(y4,((x2,y4),x3))=1,
\\
(x1,((x2,y4),x4))=(x2,((x2,y4),x4))=(x3,((x2,y4),x4))=(x4,((x2,y4),x4))=
\\
(y1,((x2,y4),x4))=(y2,((x2,y4),x4))=(y3,((x2,y4),x4))=(y4,((x2,y4),x4))=1,
\\
(x1,((x2,y4),y2))=(x2,((x2,y4),y2))=(x3,((x2,y4),y2))=(x4,((x2,y4),y2))=
\\
(y1,((x2,y4),y2))=(y2,((x2,y4),y2))=(y3,((x2,y4),y2))=(y4,((x2,y4),y2))=1,
\\
(x1,((x2,y4),y3))=(x2,((x2,y4),y3))=(x3,((x2,y4),y3))=(x4,((x2,y4),y3))=
\\
(y1,((x2,y4),y3))=(y2,((x2,y4),y3))=(y3,((x2,y4),y3))=(y4,((x2,y4),y3))=1,
\\
(x1,((x2,y4),y1))=(x2,((x2,y4),y1))=(x3,((x2,y4),y1))=(x4,((x2,y4),y1))=
\\
(y1,((x2,y4),y1))=(y2,((x2,y4),y1))=(y3,((x2,y4),y1))=(y4,((x2,y4),y1))=1,
\\
(x1,((x3,y1),x1))=(x2,((x3,y1),x1))=(x3,((x3,y1),x1))=(x4,((x3,y1),x1))=
\\
(y1,((x3,y1),x1))=(y2,((x3,y1),x1))=(y3,((x3,y1),x1))=(y4,((x3,y1),x1))=1,
\\
(x1,((x3,y1),x2))=(x2,((x3,y1),x2))=(x3,((x3,y1),x2))=(x4,((x3,y1),x2))=
\\
(y1,((x3,y1),x2))=(y2,((x3,y1),x2))=(y3,((x3,y1),x2))=(y4,((x3,y1),x2))=1,
\\
(x1,((x3,y1),x4))=(x2,((x3,y1),x4))=(x3,((x3,y1),x4))=(x4,((x3,y1),x4))=
\\
(y1,((x3,y1),x4))=(y2,((x3,y1),x4))=(y3,((x3,y1),x4))=(y4,((x3,y1),x4))=1,
\\
(x1,((x3,y1),y2))=(x2,((x3,y1),y2))=(x3,((x3,y1),y2))=(x4,((x3,y1),y2))=
\\
(y1,((x3,y1),y2))=(y2,((x3,y1),y2))=(y3,((x3,y1),y2))=(y4,((x3,y1),y2))=1,
\\
(x1,((x3,y1),y3))=(x2,((x3,y1),y3))=(x3,((x3,y1),y3))=(x4,((x3,y1),y3))=
\\
(y1,((x3,y1),y3))=(y2,((x3,y1),y3))=(y3,((x3,y1),y3))=(y4,((x3,y1),y3))=1,
\\
(x1,((x3,y1),y4))=(x2,((x3,y1),y4))=(x3,((x3,y1),y4))=(x4,((x3,y1),y4))=
\\
(y1,((x3,y1),y4))=(y2,((x3,y1),y4))=(y3,((x3,y1),y4))=(y4,((x3,y1),y4))=1,
\\
(x1,((x3,y2),x1))=(x2,((x3,y2),x1))=(x3,((x3,y2),x1))=(x4,((x3,y2),x1))=
\\
(y1,((x3,y2),x1))=(y2,((x3,y2),x1))=(y3,((x3,y2),x1))=(y4,((x3,y2),x1))=1,
\\
(x1,((x3,y2),x2))=(x2,((x3,y2),x2))=(x3,((x3,y2),x2))=(x4,((x3,y2),x2))=
\\
(y1,((x3,y2),x2))=(y2,((x3,y2),x2))=(y3,((x3,y2),x2))=(y4,((x3,y2),x2))=1,
\\
(x1,((x3,y2),x4))=(x2,((x3,y2),x4))=(x3,((x3,y2),x4))=(x4,((x3,y2),x4))=
\\
(y1,((x3,y2),x4))=(y2,((x3,y2),x4))=(y3,((x3,y2),x4))=(y4,((x3,y2),x4))=1,
\\
(x1,((x3,y2),y1))=(x2,((x3,y2),y1))=(x3,((x3,y2),y1))=(x4,((x3,y2),y1))=
\\
(y1,((x3,y2),y1))=(y2,((x3,y2),y1))=(y3,((x3,y2),y1))=(y4,((x3,y2),y1))=1,
\\
(x1,((x3,y2),y3))=(x2,((x3,y2),y3))=(x3,((x3,y2),y3))=(x4,((x3,y2),y3))=
\\
(y1,((x3,y2),y3))=(y2,((x3,y2),y3))=(y3,((x3,y2),y3))=(y4,((x3,y2),y3))=1,
\\
(x1,((x3,y2),y4))=(x2,((x3,y2),y4))=(x3,((x3,y2),y4))=(x4,((x3,y2),y4))=
\\
(y1,((x3,y2),y4))=(y2,((x3,y2),y4))=(y3,((x3,y2),y4))=(y4,((x3,y2),y4))=1,
\\
(x1,((x3,y3),x1))=(x2,((x3,y3),x1))=(x3,((x3,y3),x1))=(x4,((x3,y3),x1))=
\\
(y1,((x3,y3),x1))=(y2,((x3,y3),x1))=(y3,((x3,y3),x1))=(y4,((x3,y3),x1))=1,
\\
(x1,((x3,y3),x2))=(x2,((x3,y3),x2))=(x3,((x3,y3),x2))=(x4,((x3,y3),x2))=
\\
(y1,((x3,y3),x2))=(y2,((x3,y3),x2))=(y3,((x3,y3),x2))=(y4,((x3,y3),x2))=1,
\\
(x1,((x3,y3),x4))=(x2,((x3,y3),x4))=(x3,((x3,y3),x4))=(x4,((x3,y3),x4))=
\\
(y1,((x3,y3),x4))=(y2,((x3,y3),x4))=(y3,((x3,y3),x4))=(y4,((x3,y3),x4))=1,
\\
(x1,((x3,y3),y2))=(x2,((x3,y3),y2))=(x3,((x3,y3),y2))=(x4,((x3,y3),y2))=
\\
(y1,((x3,y3),y2))=(y2,((x3,y3),y2))=(y3,((x3,y3),y2))=(y4,((x3,y3),y2))=1,
\\
(x1,((x3,y3),y1))=(x2,((x3,y3),y1))=(x3,((x3,y3),y1))=(x4,((x3,y3),y1))=
\\
(y1,((x3,y3),y1))=(y2,((x3,y3),y1))=(y3,((x3,y3),y1))=(y4,((x3,y3),y1))=1,
\\
(x1,((x3,y3),y4))=(x2,((x3,y3),y4))=(x3,((x3,y3),y4))=(x4,((x3,y3),y4))=
\\
(y1,((x3,y3),y4))=(y2,((x3,y3),y4))=(y3,((x3,y3),y4))=(y4,((x3,y3),y4))=1,
\\
(x1,((x3,y4),x1))=(x2,((x3,y4),x1))=(x3,((x3,y4),x1))=(x4,((x3,y4),x1))=
\\
(y1,((x3,y4),x1))=(y2,((x3,y4),x1))=(y3,((x3,y4),x1))=(y4,((x3,y4),x1))=1,
\\
(x1,((x3,y4),x2))=(x2,((x3,y4),x2))=(x3,((x3,y4),x2))=(x4,((x3,y4),x2))=
\\
(y1,((x3,y4),x2))=(y2,((x3,y4),x2))=(y3,((x3,y4),x2))=(y4,((x3,y4),x2))=1,
\\
(x1,((x3,y4),x4))=(x2,((x3,y4),x4))=(x3,((x3,y4),x4))=(x4,((x3,y4),x4))=
\\
(y1,((x3,y4),x4))=(y2,((x3,y4),x4))=(y3,((x3,y4),x4))=(y4,((x3,y4),x4))=1,
\\
(x1,((x3,y4),y2))=(x2,((x3,y4),y2))=(x3,((x3,y4),y2))=(x4,((x3,y4),y2))=
\\
(y1,((x3,y4),y2))=(y2,((x3,y4),y2))=(y3,((x3,y4),y2))=(y4,((x3,y4),y2))=1,
\\
(x1,((x3,y4),y3))=(x2,((x3,y4),y3))=(x3,((x3,y4),y3))=(x4,((x3,y4),y3))=
\\
(y1,((x3,y4),y3))=(y2,((x3,y4),y3))=(y3,((x3,y4),y3))=(y4,((x3,y4),y3))=1,
\\
(x1,((x3,y4),y1))=(x2,((x3,y4),y1))=(x3,((x3,y4),y1))=(x4,((x3,y4),y1))=
\\
(y1,((x3,y4),y1))=(y2,((x3,y4),y1))=(y3,((x3,y4),y1))=(y4,((x3,y4),y1))=1,
\\
(x1,((x4,y1),x1))=(x2,((x4,y1),x1))=(x3,((x4,y1),x1))=(x4,((x4,y1),x1))=
\\
(y1,((x4,y1),x1))=(y2,((x4,y1),x1))=(y3,((x4,y1),x1))=(y4,((x4,y1),x1))=1,
\\
(x1,((x4,y1),x3))=(x2,((x4,y1),x3))=(x3,((x4,y1),x3))=(x4,((x4,y1),x3))=
\\
(y1,((x4,y1),x3))=(y2,((x4,y1),x3))=(y3,((x4,y1),x3))=(y4,((x4,y1),x3))=1,
\\
(x1,((x4,y1),x2))=(x2,((x4,y1),x2))=(x3,((x4,y1),x2))=(x4,((x4,y1),x2))=
\\
(y1,((x4,y1),x2))=(y2,((x4,y1),x2))=(y3,((x4,y1),x2))=(y4,((x4,y1),x2))=1,
\\
(x1,((x4,y1),y2))=(x2,((x4,y1),y2))=(x3,((x4,y1),y2))=(x4,((x4,y1),y2))=
\\
(y1,((x4,y1),y2))=(y2,((x4,y1),y2))=(y3,((x4,y1),y2))=(y4,((x4,y1),y2))=1,
\\
(x1,((x4,y1),y3))=(x2,((x4,y1),y3))=(x3,((x4,y1),y3))=(x4,((x4,y1),y3))=
\\
(y1,((x4,y1),y3))=(y2,((x4,y1),y3))=(y3,((x4,y1),y3))=(y4,((x4,y1),y3))=1,
\\
(x1,((x4,y1),y4))=(x2,((x4,y1),y4))=(x3,((x4,y1),y4))=(x4,((x4,y1),y4))=
\\
(y1,((x4,y1),y4))=(y2,((x4,y1),y4))=(y3,((x4,y1),y4))=(y4,((x4,y1),y4))=1,
\\
(x1,((x4,y2),x1))=(x2,((x4,y2),x1))=(x3,((x4,y2),x1))=(x4,((x4,y2),x1))=
\\
(y1,((x4,y2),x1))=(y2,((x4,y2),x1))=(y3,((x4,y2),x1))=(y4,((x4,y2),x1))=1,
\\
(x1,((x4,y2),x3))=(x2,((x4,y2),x3))=(x3,((x4,y2),x3))=(x4,((x4,y2),x3))=
\\
(y1,((x4,y2),x3))=(y2,((x4,y2),x3))=(y3,((x4,y2),x3))=(y4,((x4,y2),x3))=1,
\\
(x1,((x4,y2),x2))=(x2,((x4,y2),x2))=(x3,((x4,y2),x2))=(x4,((x4,y2),x2))=
\\
(y1,((x4,y2),x2))=(y2,((x4,y2),x2))=(y3,((x4,y2),x2))=(y4,((x4,y2),x2))=1,
\\
(x1,((x4,y2),y1))=(x4,((x4,y2),y1))=(x3,((x4,y2),y1))=(x4,((x4,y2),y1))=
\\
(y1,((x4,y2),y1))=(y2,((x4,y2),y1))=(y3,((x4,y2),y1))=(y4,((x4,y2),y1))=1,
\\
(x1,((x4,y2),y3))=(x4,((x4,y2),y3))=(x3,((x4,y2),y3))=(x4,((x4,y2),y3))=
\\
(y1,((x4,y2),y3))=(y2,((x4,y2),y3))=(y3,((x4,y2),y3))=(y4,((x4,y2),y3))=1,
\\
(x1,((x4,y2),y4))=(x2,((x4,y2),y4))=(x3,((x4,y2),y4))=(x4,((x4,y2),y4))=
\\
(y1,((x4,y2),y4))=(y2,((x4,y2),y4))=(y3,((x4,y2),y4))=(y4,((x4,y2),y4))=1,
\\
(x1,((x4,y3),x1))=(x2,((x4,y3),x1))=(x3,((x4,y3),x1))=(x4,((x4,y3),x1))=
\\
(y1,((x4,y3),x1))=(y2,((x4,y3),x1))=(y3,((x4,y3),x1))=(y4,((x4,y3),x1))=1,
\\
(x1,((x4,y3),x3))=(x2,((x4,y3),x3))=(x3,((x4,y3),x3))=(x4,((x4,y3),x3))=
\\
(y1,((x4,y3),x3))=(y2,((x4,y3),x3))=(y3,((x4,y3),x3))=(y4,((x4,y3),x3))=1,
\\
(x1,((x4,y3),x2))=(x2,((x4,y3),x2))=(x3,((x4,y3),x2))=(x4,((x4,y3),x2))=
\\
(y1,((x4,y3),x2))=(y2,((x4,y3),x2))=(y3,((x4,y3),x2))=(y4,((x4,y3),x2))=1,
\\
(x1,((x4,y3),y2))=(x2,((x4,y3),y2))=(x3,((x4,y3),y2))=(x4,((x4,y3),y2))=
\\
(y1,((x4,y3),y2))=(y2,((x4,y3),y2))=(y3,((x4,y3),y2))=(y4,((x4,y3),y2))=1,
\\
(x1,((x4,y3),y1))=(x2,((x4,y3),y1))=(x3,((x4,y3),y1))=(x4,((x4,y3),y1))=
\\
(y1,((x4,y3),y1))=(y2,((x4,y3),y1))=(y3,((x4,y3),y1))=(y4,((x4,y3),y1))=1,
\\
(x1,((x4,y3),y4))=(x2,((x4,y3),y4))=(x3,((x4,y3),y4))=(x4,((x4,y3),y4))=
\\
(y1,((x4,y3),y4))=(y2,((x4,y3),y4))=(y3,((x4,y3),y4))=(y4,((x4,y3),y4))=1,
\\
(x1,((x4,y4),x1))=(x2,((x4,y4),x1))=(x3,((x4,y4),x1))=(x4,((x4,y4),x1))=
\\
(y1,((x4,y4),x1))=(y2,((x4,y4),x1))=(y3,((x4,y4),x1))=(y4,((x4,y4),x1))=1,
\\
(x1,((x4,y4),x3))=(x2,((x4,y4),x3))=(x3,((x4,y4),x3))=(x4,((x4,y4),x3))=
\\
(y1,((x4,y4),x3))=(y2,((x4,y4),x3))=(y3,((x4,y4),x3))=(y4,((x4,y4),x3))=1,
\\
(x1,((x4,y4),x2))=(x2,((x4,y4),x2))=(x3,((x4,y4),x2))=(x4,((x4,y4),x2))=
\\
(y1,((x4,y4),x2))=(y2,((x4,y4),x2))=(y3,((x4,y4),x2))=(y4,((x4,y4),x2))=1,
\\
(x1,((x4,y4),y2))=(x2,((x4,y4),y2))=(x3,((x4,y4),y2))=(x4,((x4,y4),y2))=
\\
(y1,((x4,y4),y2))=(y2,((x4,y4),y2))=(y3,((x4,y4),y2))=(y4,((x4,y4),y2))=1,
\\
(x1,((x4,y4),y3))=(x2,((x4,y4),y3))=(x3,((x4,y4),y3))=(x4,((x4,y4),y3))=
\\
(y1,((x4,y4),y3))=(y2,((x4,y4),y3))=(y3,((x4,y4),y3))=(y4,((x4,y4),y3))=1,
\\
(x1,((x4,y4),y1))=(x2,((x4,y4),y1))=(x3,((x4,y4),y1))=(x4,((x4,y4),y1))=
\\
(y1,((x4,y4),y1))=(y2,((x4,y4),y1))=(y3,((x4,y4),y1))=(y4,((x4,y4),y1))=1,
\\
((x1*x2*x3*x4,y1),x1)*((y1,x1),y2)*((x1,y2),x3)*((y2,x3),y4)*((x3,y4),x2)*
\\
((y4,x2),y3)=((x2,y3),x1*x2*x3*x4)*((y3,x1*x2*x3*x4),y1),
\\
((x4,y2),x1*x2*x3*x4)*((y2,x1*x2*x3*x4),y3)*((x1*x2*x3*x4,y3),x2)*
\\
((y3,x2),y1*y2*y3*y4)*((x2,y1*y2*y3*y4),x1)=
\\
((y1*y2*y3*y4,x1),y4)*((x1,y4),x4)*((y4,x4),y2)
>;}
\smallskip

\f {Construct the largest 2-quotient group of $P$ having lower exponent-$2$ class at most 200 as group in the category GrpPC:}
\smallskip
\\
{\tt P,q:=pQuotient(G,2,200);}\smallskip

\f{Order of $\P$ (The result shows that $|P|=2^{59}$, and so $\P\cong P$ has order $2^{59}$):}\smallskip
\\
{\tt FactoredOrder(P);}\smallskip
\\
\f{Construct the subgroup $\H$:}\smallskip
\\
{\tt x1:=x1@q; x2:=x2@q; x3:=x3@q; x4:=x4@q;
\\
y1:=y1@q; y2:=y2@q; y3:=y3@q; y4:=y4@q; r:=r@q;
\\
H:=sub<P|x1,x2,x3,x4,y1,y2,y3,y4>;
\\
FactoredOrder(H);
\\
IsNormal(P,H);
}
\\
\f{$r$ is an automorphism of $\H$ of order $8$:}\smallskip
\\
{\tt
Order(r);
\\
$\sharp$(Centraliser(P,H) meet sub<P|r>);
}\smallskip
\\
\f{$\H$ is a $4$-dimensional mixed dihedral group relative to $X=\lg x_1,x_2,x_3,x_4\rg$ and $Y=\lg y_1,y_2,y_3,y_4\rg$:}\smallskip
\\
{\tt GroupName(H/DerivedSubgroup(H));
\\
GroupName(sub<P|x1,x2,x3,x4>);
\\
GroupName(sub<P|y1,y2,y3,y4>);}\smallskip
\\
\f{$\a_1,\a_2,\b_1,\b_2$ are automorphisms of $\H$:}\smallskip
\\
{\tt
alpha1:=hom<H->H|x1->x1*x2, x2->x2*x3, x3->x3*x4, x4->x1*x2*x3,
\\
y1->y1, y2->y2, y3->y3, y4->y4>;
\\
$\sharp$Kernel(alpha1);Image(alpha1) eq H;
\\
alpha2:=hom<H->H|x1->x1, x2->x2, x3->x3, x4->x4,
\\
y1->y1*y2, y2->y2*y3, y3->y3*y4, y4->y1*y2*y3>;
\\
$\sharp$Kernel(alpha2);Image(alpha2) eq H;
\\
beta1:=hom<H->H|x1->x2, x2->x3, x3->x4, x4->x1*x2*x3*x4,
\\
y1->y1, y2->y2, y3->y3, y4->y4>;
\\
$\sharp$Kernel(beta1);Image(beta1) eq H;
\\
beta2:=hom<H->H|x1->x1, x2->x2, x3->x3, x4->x4,
\\
y1->y2, y2->y3, y3->y4, y4->y1*y2*y3*y4>;
\\
$\sharp$Kernel(beta1);Image(beta1) eq H;
}\smallskip
\\
\f{Neither $\s_1$ nor $(r^4)^X$ is an automorphism of $\H$:}\smallskip
\\
{\tt
sigma1:=hom<H->H | x1->x1, x2->x1*x3, x3->x2*x3, x4->x2*x4, y1->y1, y2->y2, y3->y3, y4->y4>;
\\
$(r^4)^X$:=hom<H->H | x1->x4, x2->x3, x3->x2, x4->x1, y1->y1, y2->y2, y3->y3, y4->y4>;
}\smallskip
\\
\f{Stabiliser of $X$ in $\P$:}\smallskip
\\
{\tt
Stab:=sub<P|x1,x2,x3,x4,r\textasciicircum2>;
\\
FactoredOrder(Stab);
}\smallskip
\\
\f{Transitive subgroups of $\P$ of order $2^{58}$:}\smallskip
\\
{\tt
S:=MaximalSubgroups(P);
\\
T:=\{\};
\\
for i in \{1..$\sharp$S\} do
\\
if $\sharp$(S[i]`subgroup meet Stab) eq 2\textasciicircum5 then
\\
Include($^\sim$T,S[i]`subgroup);
\\
end if;
\\
end for;
$\sharp$T;
}\smallskip
\\
\f{Transitive subgroups of $\P$ of order $2^{57}$:}\smallskip
\\
{\tt
MT:=\{\};
\\
for x in T do
\\
Mx:=MaximalSubgroups(x);
\\
for i in \{1..$\sharp$Mx\} do
\\
if $\sharp$(Mx[i]`subgroup meet Stab) eq 2\textasciicircum4 then
\\
Include($^\sim$MT,Mx[i]`subgroup);
\\
end if;
\\
end for;
\\
end for;
\\
$\sharp$MT;
}
\smallskip
\\
\f{Transitive subgroups of $\P$ of order $2^{56}$:}\smallskip
\\
{\tt
MTT:=\{\};
\\
for x in MT do
\\
Mx:=MaximalSubgroups(x);
\\
for i in \{1..$\sharp$Mx\} do
\\
if $\sharp$(Mx[i]`subgroup meet Stab) eq 2\textasciicircum3 then
\\
Include($^\sim$MTT,Mx[i]`subgroup);
\\
end if;
\\
end for;
\\
end for;
\\
$\sharp$MTT;
}
\smallskip
\\
\f{Transitive subgroups of $\P$ of order $2^{55}$:}\smallskip
\\
{\tt
MTTT:=\{\};
\\
for x in MTT do
\\
Mx:=MaximalSubgroups(x);
\\
for i in \{1..$\sharp$Mx\} do
\\
if $\sharp$(Mx[i]`subgroup meet Stab) eq 2\textasciicircum2 then
\\
Include($^\sim$MTTT,Mx[i]`subgroup);
\\
end if;
\\
end for;
\\
end for;
\\
$\sharp$MTTT;
}
\smallskip
\\
\f{Transitive subgroups of $\P$ of order $2^{54}$:}\smallskip
\\
{\tt
MTTTT:=\{\};
\\
for x in MTTT do
\\
Mx:=MaximalSubgroups(x);
\\
for i in \{1..$\sharp$Mx\} do
\\
if $\sharp$(Mx[i]`subgroup meet Stab) eq 2 then
\\
Include($^\sim$MTTTT,Mx[i]`subgroup);
\\
end if;
\\
end for;
\\
end for;
\\
$\sharp$MTTTT;
}
\smallskip
\\
\f{No regular subgroups of $\P$:}\smallskip
\\
{\tt
MTTTTT:=\{\};
\\
for x in MTTTT do
\\
Mx:=MaximalSubgroups(x);
\\
for i in \{1..$\sharp$Mx\} do
\\
if $\sharp$(Mx[i]`subgroup meet Stab) eq 1 then
\\
Include($^\sim$MTTTTT,Mx[i]`subgroup);
\\
end if;
\\
end for;
\\
end for;
\\
$\sharp$MTTTTT;
}
\bigskip\bigskip

\section*{Acknowledgements}

The first author has been supported by the Croatian
Science Foundation under the project 6732. The second author is grateful for Australian Research Council Discovery Project Grant DP230101268. The third author was supported by the National Natural Science Foundation of China (12071023, 12161141005) and the 111 Project of China (B16002).

\end{document}